\documentclass{amsart}

\newtheorem{theorem}{Theorem}
\newcommand{\bt}{\begin{theorem}}
\newcommand{\et}{\end{theorem}}

\newtheorem{lemma}{Lemma}
\newcommand{\bl}{\begin{lemma}}
\newcommand{\el}{\end{lemma}}

\newtheorem{problem}{Problem}
\newcommand{\bprob}{\begin{problem}}
\newcommand{\eprob}{\end{problem}}

\newcommand{\Q}{\ensuremath{ \mathbf{Q} }}

\usepackage{amssymb,latexsym}
\usepackage[all]{xy}

\begin{document}
 
\title{Geometric progressions in syndetic sets}
\author{Melvyn B. Nathanson}
\address{Department of Mathematics\\Lehman College (CUNY)\\Bronx, NY 10468}
\email{melvyn.nathanson@lehman.cuny.edu}

\thanks{This research was supported in part by the International Centre for Theoretical Sciences (ICTS) 
 during a visit for the program -  Workshop on Additive Combinatorics (Code:ICTS/Prog-wac2020/02). }


  
 \date{\today}

 \begin{abstract}
An open problem about finite geometric progressions in syndetic sets leads to a family of diophantine 
equations related to the commutativity of translation and multiplication by squares.   
\end{abstract}

 \maketitle

\section{The problem of geometric progressions in syndetic sets}
The set  $A = \{ a_i : i = 1,2,3,\ldots \}$ of positive integers is \emph{strictly increasing} 
if  $a_1 < a_2 < a_3 < \cdots$. 
The set $A$ has \emph{gaps bounded by $\ell$} if $a_{i+1}-a_i \leq \ell$ for all $i= 1,2,3,\ldots$.
The set $A$ has \emph{bounded gaps} if it has gaps bounded by $\ell$ for some $\ell$.
A set with bounded gaps is also called \emph{syndetic}. 
Beiglb{\"o}ck, Bergelson, Hindman,  and Strauss~\cite{beig-berg-hind-stra06} 
asked if every syndetic set of integers contains arbitrarily long finite geometric progressions.  
For recent progress on this still unsolved problem, 
see Glasscock,  Koutsogiannis, and Richter~\cite{glas19} and Patil~\cite{pati19}.
For sets not containing finite geometric progressions, 
see McNew~\cite{mcne20} and 
Nathanson and O'Bryant~\cite{nath2013-149,nath2014-154,nath2015-159}. 

It is not even known if a syndetic set contains a three-term geometric progression, that is, 
a subset  of the form $\{a,ax,ax^2\}$ for integers $a$ and $x$ with $x \geq 2$. 
Indeed, it is not known if a syndetic set must contain a subset of the form $\{a,ax^2\}$ 
for some integer  $x \geq 2$. 
Patil~\cite{pati19} recently proved the following special case: 
Every strictly increasing set of positive integers with gaps bounded by 2  
contains infinitely many subsets $\{a,ax_a^2\}$ with $x_a \geq 2$.  
His proof uses the following beautiful polynomial  identity:
\begin{equation}                \label{geo-syn:patil}
a(4a+3)^2+1 = (a+1)(4a+1)^2.
\end{equation}
The core of the argument is the following result.

\bt[Patil]                      \label{PowerShift:theorem:a+1}
Let $A$ be a set of positive integers such that $\{a,a+1\}\subseteq A$ for infinitely many $a\in A$,
and $\{a,a+1\} \cap A \neq \emptyset$ for all positive integers $a$.
The set $A$ contains infinitely many subsets of the form $\{a,ax_a^2\}$ with integers $x_a \geq 2$.  
\et

The object of this note is to construct an infinite family of arithmetic identities 
that will generalize Theorem~\ref{PowerShift:theorem:a+1}.

\section{Translations and dilations by squares}

Let $a$, $k$,  $\ell$, $m$ and $n$  be positive integers with $m \geq 2$ and $n \geq 2$. 
Does there exist an integer $b$, or do there exist infinitely many integers $b$, such that 
\begin{enumerate}
\item[(i)] $a$ multiplied by an $m$th power and translated by $k$ is equal to $b$, and 
\item[(ii)] $a$ translated by $\ell$  and multiplied by an $n$th power is also equal to $b$?
\end{enumerate}
Thus, the problem asks about the existence of a pair of positive integers $(x,y)$, 
or of infinitely many pairs of positive integers $(x,y)$, such that 
\[
ax^m + k = (a+\ell)y^n.  
\]
A necessary condition is
\[
k \equiv 0 \pmod{\gcd(a,\ell)}.
\] 
  
We solve this problem in the case of equal translations and of dilations by squares, 
that is, $k = \ell$ and $m=n=2$.

\bt               \label{PowerShift:theorem:PowerShift}
Let $a$ and $k$ be positive integers such that 
\[
d = a(a+k) 
\]
is not a square.  Let $(u,v)$ be a positive integral solution of the Pell's equation 
\begin{equation}                            \label{PowerShift:Pell}
u^2 - dv^2 = 1.
\end{equation} 
If
\begin{equation}                            \label{PowerShift:squares-solution-x}
x= u+av+kv  \qquad \text{and} \qquad y = u+av
\end{equation}
then 
\begin{equation}                            \label{PowerShift:squares-solution-y}
ax^2 + k = (a+k)y^2.
\end{equation}
\et

\begin{proof}
This is a simple calculation.  Insert~\eqref{PowerShift:squares-solution-x} 
into~\eqref{PowerShift:squares-solution-y} and use relation~\eqref{PowerShift:Pell}.  
This completes the proof. 
\end{proof}

The solution~\eqref{PowerShift:squares-solution-x} was obtained as follows.  
A rearrangement of~\eqref{PowerShift:squares-solution-y} is
\[
(a+k)y^2 - ax^2 = k.
\]
Writing 
\[
z = (a+k)y  \qquad \text{and} \qquad d = a(a+k)
\]
we obtain 
\[
z^2 - dx^2 = \left( (a+k)y\right)^2 - a(a+k)x^2 = k(a+k).
\]
Equivalently, in the real quadratic field $\Q(\sqrt{d})$ with norm $N\left( u+v\sqrt{d}\right) = u^2 - dv^2$, 
we have  
\[
N\left( z+x\sqrt{d}\right)  = k(a+k).
\]
Consider a particular solution  
\[
N\left( z_0 + x_0\sqrt{d}\right) = k(a+k).
\]
If $(u,v)$ is a solution of Pell's equation~\eqref{PowerShift:Pell}, 
then $N\left( u+v\sqrt{d}\right) = 1$.  
We have 
\[
z_1 + x_1\sqrt{d} = \left( z_0 + x_0\sqrt{d}\right) \left(u+v\sqrt{d} \right) 
= (z_0u + x_0vd ) + (z_0v+x_0u)\sqrt{d}. 
\]
The multiplicativity of the norm gives 
\begin{align*}
N\left( z_1 + x_1\sqrt{d}\right) 
& = N\left( z_0 + x_0\sqrt{d}\right) N\left( u+v\sqrt{d}\right) \\ 
& = k(a+k).
\end{align*}
Choosing the particular solution $(x_0,y_0) = (1,1)$ 
and $(x_0,z_0) = (1,a+k)$, we obtain 
\[
x_1 = z_0v+x_0u = u + av + kv
\]
and 
\begin{align*}
(a+k)y_1 & = z_1 = z_0u + x_0vd \\
& = (a+k)u + va(a+k) \\
& =  (a+k)(u + av) 
\end{align*}
and so 
\[
y_1 = u + av.
\]
This gives~\eqref{PowerShift:squares-solution-x}.

\bl                                                               \label{PowerShift:lemma:a+k}
Let $k$ be a positive integer.
There are only finitely many positive integers $a$ such that $a(a+k)$ is a square.
\el

\begin{proof}
Let 
\[
a(a+k) = r^2.
\]
We write
\[
a = b^2 c
\]
where $c$ is square-free.  
The equation 
\[
a(a+k) =  b^2 c \left( b^2 c + k \right) = r^2
\]
implies that $b$ divides $r$.  Let $r = bs$.  We obtain 
\[
b^2 c \left( b^2 c + k \right) = r^2 = b^2 s^2
\]
and  
\[
 c \left( b^2 c + k \right) = s^2.  
\]
The integer $c$ is square-free, and so  $c$ divides $s$.  Let $s = ct$.  We obtain 
\[
 c \left( b^2 c + k \right) = s^2 = c^2 t^2
\]
and  
\[
b^2 c + k =  ct^2.
\]
Therefore, $c$ divides $b^2 c + k$ and so $c$ divides $k$.  Let $k = c\ell$.
We obtain 
\[
b^2 c + c\ell  = ct^2 
\] 
and 
\[
b^2 + \ell  = t^2.
\]
For every positive integer $\ell$, there exist only finitely many pairs of  positive integers $(b,t)$ 
such that $t^2 - b^2 = \ell$, and so, for every positive divisor $\ell$ of $k$, 
there are only 
finitely many positive integers $a$ such that $a(a+k)$ is a square.
This completes the proof. 
\end{proof}

\bt                       \label{PowerShift:theorem:a+k}
Let $k$ be a positive integer.  
Let $A$ be a set of positive integers such that $\{a,a+k\}\subseteq A$ for infinitely many $a\in A$,
and $\{a,a+k\} \cap A \neq \emptyset$ for all positive integers $a$.
The set $A$ contains infinitely many subsets of the form $\{a,ax^2\}$ for some integer $x \geq 2$.  
\et

\begin{proof}
By Lemma~\ref{PowerShift:lemma:a+k}, the integer $a(a+k)$ is a square for only finitely many 
$a\in A$, and so we can assume that $a(a+k)$ is not a square for all $a \in A$.  

Let $\{a,a+k\}\subseteq A$.   Apply~\eqref{PowerShift:squares-solution-x} in  Theorem~\ref{PowerShift:theorem:PowerShift} 
to construct positive integers $x$ and $y$.  Let $b = ax^2$.

If $b\in A$, then $A$ contains the set $\{a,ax^2\}$.  If $b = ax^2 \notin A$, then 
\[
b+k  = ax^2 + k = (a+k)y^2\in A
\]
and  $A$ contains the set $\{a+k, (a+k) y^2\}$. 
This completes the proof.  
\end{proof}

Note that Patel's theorem is the case $k=1$ of this result.  

\section{Open problems}

\bprob
Determine the triples of positive integers $(a, k,\ell)$ with $k \neq  \ell$ 
such that there exists one pair or there exist infinitely many pairs of squares $(x^2,y^2)$ 
with 
\[
ax^2 + k = (a+ \ell)y^2.
\]
\eprob

\bprob
Determine the triples of positive integers $(a, k,\ell)$ 
such that there exists one pair or there exist infinitely many pairs of cubes $(x^3,y^3)$ 
with 
\[
ax^3 + k = (a+ \ell)y^3.
\]
\eprob

\bprob
Determine the quintuples of positive integers $(a, k,\ell,m,n)$ 
such that there exists one pair or there exist infinitely many pairs of integer powers $(x^m,y^n)$ 
with 
\[
ax^m + k = (a+ \ell)y^n.
\]
\eprob

I wish to thank B.~R. Patil for introducing me to this subject at the 
ICTS Workshop on Additive Combinatorics in Bangalore in March, 2020.

Added June 27, 2020:  
I also thank Jaitra Chattopadhyay for observing that, in the first version of this paper, 
the condition ``$d=a(a+k)$ is not a square" 
had been omitted in the statement of Theorem~\ref{PowerShift:theorem:PowerShift}.

\def\cprime{$'$} \def\cprime{$'$} \def\cprime{$'$}
\providecommand{\bysame}{\leavevmode\hbox to3em{\hrulefill}\thinspace}
\providecommand{\MR}{\relax\ifhmode\unskip\space\fi MR }
\providecommand{\MRhref}[2]{%
  \href{http://www.ams.org/mathscinet-getitem?mr=#1}{#2}
}
\providecommand{\href}[2]{#2}

\end{document}